\documentclass[11pt]{amsart}
\usepackage{amsfonts}
\usepackage{amsmath}
\usepackage{amsthm}
\usepackage{bm,bbm}
\usepackage[dvips]{graphicx}
\usepackage{caption}
\usepackage[numbers]{natbib}
\usepackage{color}
\usepackage{soul}

\usepackage{bigints}

\usepackage{bm}

\usepackage{mathtools}
\usepackage{commath}

\RequirePackage[colorlinks=true, linkcolor = blue,citecolor=black,urlcolor=blue,backref=page]{hyperref}

\numberwithin{equation}{section}

\allowdisplaybreaks

\theoremstyle{plain}
\newtheorem{theorem}{Theorem}[section]
\newtheorem{lemma}[theorem]{Lemma}

\newtheorem{corollary}[theorem]{Corollary}
\theoremstyle{definition}
\newtheorem{definition}[theorem]{Definition}

\newtheorem{example}[theorem]{Example}

\def\ms{\mathsf}
\def\beqn{\begin{equation}}
\def\beqn*{$$}
\def\eeqn{\end{equation}}

\def\ms{\mathsf}

\newcommand{\bx}{{\bf x}}

\newcommand{\bz}{{\bf z}}

\def\P{\mathbb{P}}

\def\E{\mathbb{E}}
\def\Pn{\mathcal P_n}

\newcommand{\reals}{{\mathbb R}}
\newcommand{\bbr}{\reals}
\newcommand{\R}{\reals}
\newcommand{\bbn}{{\mathbb N}}

\newcommand{\vep}{\varepsilon}
\newcommand{\bbz}{\protect{\mathbb Z}}

\newcommand{\X}{{\mathcal{X}}}
\newcommand{\Y}{{\mathcal{Y}}}

\newcommand{\B}{\mathcal B}

\newcommand{\btheta}{{\bm \theta}}

\newcommand{\one}{{\mathbbm 1}}

\newcommand{\remove}[1]{}

\newcommand{\C}{\check{C}}

\def\enp{\end{proof}}
\def\bel{\begin{lemma}}
\def\bep{\begin{proof}}
\def\enl{\end{lemma}}
\newcommand{\M}{\mathcal M}

\def\T{\mathbb{T}}
\def\bbY{\mathbb{Y}}
\newcommand{\bkn}{b_{k,n}}

\begin{document}

\date{\today}
\bibliographystyle{abbrv}

\renewcommand{\baselinestretch}{1.05}

\title[Critical faces above the vanishing threshold]
{Limit theorems for critical faces above the vanishing threshold}

\author{Zifu Wei}
\address{Department of Statistics\\
Purdue University \\
West Lafayette, 47907, USA}
\email{wei296@purdue.edu}

\author{Takashi Owada}
\address{Department of Statistics\\
Purdue University \\
West Lafayette, 47907, USA}
\email{owada@purdue.edu}

\author{D.Yogeshwaran}
\address{Theoretical Statistics and Mathematics Unit\\
Indian Statistical Institute \\
Bangalore, India}
\email{d.yogesh@isibang.ac.in}

\thanks{TO's research was partially supported by the AFOSR grant FA9550-22-1-0238. DY was funded through CPDA of the Indian Statistical Institute and SERB-MATRICS: MTR-2020-000470.}

\keywords{$\M_0$-convergence, point process, Morse critical points, stochastic geometry.  \vspace{.5ex}}

\begin{abstract}
We investigate convergence of point processes associated with critical faces for a \v{C}ech filtration built over a homogeneous Poisson point process in the $d$-dimensional flat torus. The convergence of our point process is established in terms of the $\mathcal M_0$-topology, when the connecting radius of a \v{C}ech complex decays to $0$, so slowly that critical faces are even less likely to occur than those in the regime of  threshold for homological connectivity. We also obtain a series of limit theorems for positive and negative critical faces, all of which are considerably analogous to those for critical faces. 
\end{abstract}

\maketitle

%%%%%%%%%%%%%%%%%%%%%%%%%%%%%%%%%%%%%%%%%%%%%%
%%%% Main text entry area:

\section{Introduction}

The main objective of this paper is to study geometric and topological configuration of point processes associated to the spatial distribution of critical faces. We will introduce  necessary notions informally and outline our results. More rigorous definitions and statements are postponed to Section \ref{sec:setup.results}. 

Consider the $d$-dimensional flat torus $\T^d =\R^d/\bbz^d$ (with a periodic boundary), that is equipped with the toroidal metric 
$$
\ms{dist}(x,y) = \min_{z\in \bbz^d}\| x-y+z \|, \ \ \ x,y \in \T^d,
$$
where $\|\cdot\|$ denotes the Euclidean norm. Given a finite point set $\X\subset \R^d$, the {\em distance function} $d_\X:\R^d\to[0,\infty)$, that will be examined in this paper, is defined by                                                        
\begin{equation}  \label{e:min.type.dist.func}
d_\X(x)=\min_{y\in \X}\ms{dist}(x,y),  \ \ x\in \R^d. 
\end{equation}
Morse theory informs us that the homology of the {\em sub-level set} $C(\X,r) = \{x \in \R^d : d_\X(x) \leq r\}$ changes at {\em critical levels} of the distance function.  The distance function is not smooth, but it is of min-type and so critical levels and critical points can be defined; see \cite{gershkovich:rubinstein:1997}.  In particular,  a critical point of index $k\in \{ 1,\ldots,d\}$ can be generated by a $(k+1)$-tuple $\Y \subset \X$. Such a critical point shall be denoted by $c(\Y)$ and a critical value or level by $\rho(\Y)$; see Definition \ref{def:crit.point}.  

This connection with Morse theory has been exploited  to study thresholds for vanishing of homology of $C(\Pn,r)$ (see  \cite{bobrowski:adler:2014,bobrowski:mukherjee:2015,bobrowski:weinberger:2017,bobrowski:oliveira:2019,bobrowski:2022}), where $\Pn$ is a homogeneous Poisson point process with intensity $n$ in $\T^d$. To better understand the behaviour of homology above the vanishing threshold, we investigate point processes of the corresponding critical levels and points. 

We now introduce some notation to outline our results. Since we are working with the torus, we do not consider sub-level sets for large levels and hence we focus only on $(C(\Pn,r))_{r\in [0,R_n]}$, where $R_n$ is a sequence of positive numbers tending to $0$ slowly enough. Given such $R_n$, we are interested in the stochastic features of the following point process induced by critical points and values:
\begin{equation}  \label{e:def.eta.kn.intro}
\eta_{k,n}:= \sum_{\Y\subset \Pn, \, |\Y|=k+1} s_{k,n}(\Y, \Pn)\,\delta_{(c(\Y), \,n\omega_d \rho(\Y)^d-a_n)}, 
\end{equation}
where $\delta_{(x_1,x_2)}$ is the Dirac measure at $(x_1,x_2)\in \T^d\times (-\infty,\infty]$, and $s_{k,n}(\Y,\Pn)$ is an indicator function, requiring that $\Y$ forms a critical $k$-face, such that $\rho(\Y)\le R_n$. Furthermore, $\omega_d$ denotes volume of the $d$-dimensional unit ball, and $a_n$ is a properly defined centering sequence. Closely related to the above point process is the statistics of total number of critical $k$-faces with large critical values. Given  $u\in \R$, let $(r_n(u))_{n\ge1}$ be a sequence of positive numbers defined by 
\begin{equation}  \label{e:def.rn.u}
r_n(u) = \Big(  \frac{a_n+u}{n\omega_d}\Big)^{1/d},  \ \ \ n\ge1. 
\end{equation}
Fixing $u=u_0$ and assuming  $r_n(u_0)\to0$, $r_n(u_0)/R_n\to0$ as $n\to\infty$, we also study the behavior of the statistics 
\begin{equation}  \label{e:def.Gkn.intro}
G_{k,n} := \sum_{\Y\subset \Pn, \, |\Y|=k+1} s_{k,n}(\Y, \Pn)\,\one \big\{ \rho(\Y)\ge r_n (u_0)\big\}, 
\end{equation}
which counts the number of critical $k$-faces whose critical values are between $r_n(u_0)$ and $R_n$. Note that $G_{k,n} = \eta_{k,n}(\T^d \times [u_0,\infty))$.

In this paper, we put our focus on the dense regime i.e., $nr_n(u)^d \to \infty$. Poisson process approximation for the process $\eta_{k,n}$ in the regime
\begin{equation}
 \label{e:def.an.Poisson.intro}   
 a_n = \log n +(k-1)\log \log n + \ms{const},
\end{equation}
as well as  Poisson convergence for $G_{k,n}$ in \eqref{e:def.Gkn.intro}, were established in Theorem 6.1 of \cite{bobrowski:schulte:yogeshwaran:2022} and Theorem 8.1 of \cite{bobrowski:2022},  respectively. These results were crucial to describing the behaviour of homology at the vanishing threshold.

We place ourselves above the Poisson regime \eqref{e:def.an.Poisson.intro}, or equivalently, above the vanishing threshold. More accurately, we are interested in the asymptotics of $G_{k,n}$, in the case that 
$$
r_n(u) \to 0, \ \ \ r_n(u) \gg \Big( \frac{\log n +(k-1)\log \log n}{n\omega_d} \Big)^{1/d}. 
$$
Equivalently, we consider  the following  centering term $(a_n)$ in \eqref{e:def.eta.kn.intro}:
%$u = u(n) \to \infty$ and asymptotics for $\eta_{k,n}$ with $(a_n)$ as follows.
\begin{equation}  \label{e:def.an.M0}
a_n-\log n -(k-1)\log \log n\to\infty, \ \ \ a_n=o(n), \ \ \ n\to\infty. 
\end{equation}
Under this assumption, the process $\eta_{k,n}$ converges to the null measure (i.e., the measure that assigns zeros to all Borel measurable sets in $\T^d\times (-\infty,\infty]$), and $G_{k,n}$ converges to a degenerate zero random variable (\cite{bobrowski:weinberger:2017}). Our main results - Theorems \ref{t:M0.eta.kn} and \ref{t:vague.eta.kn} below - quantify the rate of convergence to the null measure and a zero random variable using appropriate notions. More specifically, we identify a sequence {$b_{k,n} := na_n^{k-1}e^{-a_n} \to 0$} to show that 
\begin{equation}  \label{e:eta.kn.prob.measure.intro}
\big( b_{k,n}^{-1}\P(\eta_{k,n}\in \cdot), \, n\ge1 \big) %\ \ \text{ and } \ \ 
\end{equation}
and 
\begin{equation}   \label{e:G.kn.prob.measure.intro}
\big( b_{k,n}^{-1}\P(G_{k,n}\in \cdot), \, n\ge1 \big) 
\end{equation}
tend to non-degenerate limits as $n\to\infty$, while identifying the limits themselves explicitly. Up to a constant (dependent on $k$), the limit in \eqref{e:eta.kn.prob.measure.intro} is supported on singletons with a  product density independent of $k$. The topology underlying the convergence of \eqref{e:eta.kn.prob.measure.intro} is called \emph{$\M_0$-topology}, a useful notion  for dealing with convergence of probability measures defined on a complete and separable metric space. In recent times, the notion of $\M_0$-topology has been used for the study of geometrically and/or topologically rare events (\cite{owada:2022, hirsch:kang:owada:2023, owada:2023}), as well as regular variation of point processes and stochastic processes (\cite{hult:lindskog:2006a, hult:samorodnitsky:2010, fasen:roy:2016, segers:zhao:meinguet:2017}). For further analyses, we recall that critical faces can be divided into positive and negative critical faces (see \cite{bobrowski:2022}). Loosely speaking, positive critical $k$-faces will create a (nontrivial) $k$-dimensional cycle in the $k$th homology group of  $C(\X,r)$, while negative critical $k$-faces terminate a $(k-1)$-dimensional cycle. According to Propositions 4.2 and 5.1 in \cite{bobrowski:2022}, for each $1\le k \le d-2$, the vanishing thresholds for critical $k$-faces, positive critical $k$-faces, and negative critical $(k+1)$-faces all coincide with one another (with high probability). Extending this further, we can prove  that even when the radius $(r_n(u))_{n\ge1}$ satisfies condition  \eqref{e:def.an.M0}, the functionals of positive critical $k$-faces and  negative critical $(k+1)$-faces  exhibit  asymptotic results similar to those for critical $k$-faces; see Corollary \ref{cor:pos.neg.faces} for more details. 
The Poisson convergence results for positive and negative critical faces in the regime \eqref{e:def.an.Poisson.intro} were used to understand the vanishing threshold for homology in \cite{bobrowski:2022}. Similarly, with the equivalence between positive (resp.~negative) critical faces and birth (resp. death) times in persistence diagrams, our results %for $\eta^+_{k,n}$ (resp. $\eta^{-}_{k,n}$) 
can be rephrased in terms of the birth  and  death times within the interval $[r_n(u),R_n]$. In other words, we can quantify the distribution of 'noisy barcodes' in the persistence diagram above the homological connectivity regime. 

As a final remark, we point out that in addition to these studies on the vanishing homology, together with the related Poisson convergence, there have also been a number of attempts at deducing other types of limit theorems, including the central limit theorems (\cite{bobrowski:adler:2014, bobrowski:mukherjee:2015, yogeshwaran:adler:2015}) and the large deviation principle in \cite{hirsch:owada:2023}, when the radius $(r_n(u))_{n\ge1}$ belongs to the sparse regime (i.e., $nr_n(u)^d\to0$) or critical regime (i.e., $nr_n(u)^d\to c\in (0,\infty)$). 

The remainder of the paper is structured as follows. Section \ref{sec:setup.results} gives a precise setup for the processes $\eta_{k,n}$ and $G_{k,n}$. After that, the convergence results \eqref{e:eta.kn.prob.measure.intro} and \eqref{e:G.kn.prob.measure.intro} will be formally described in Theorems \ref{t:M0.eta.kn} and \ref{t:vague.eta.kn}. {Corollary \ref{cor:pos.neg.faces} states the convergence results for positive and negative critical faces.} All the proofs are postponed to Section \ref{sec:proofs}. The main machinery in our proof is given by Theorems 4.1 and 6.1 of  \cite{bobrowski:schulte:yogeshwaran:2022}; this will help us to show that under the $\M_0$-topology, the process $\eta_{k,n}$ can be approximated by some Poisson point process whose intensity measure tends to the null measure as $n\to\infty$. We then rely upon the estimates in \cite{bobrowski:2022} to approximate the spatial distributions of positive and negative critical faces by those of critical faces. 

\section{Setup and main results}  \label{sec:setup.results}

To better understand critical faces, it is important to define a \emph{\v{C}ech complex}, which is one of the most studied geometric complexes (see \cite{ghrist:2014}).

\begin{definition}  \label{def:cech}
Given a point set $\X=\{ x_1,\dots,x_n \}\subset \T^d$ and $r>0$, the \v{C}ech complex $\C(\X,r)$ is defined as follows:
\begin{itemize}
\item The $0$-simplices are the points in $\X$. 
\item For each $m\ge1$, $\{x_{i_0},\dots,x_{i_m}\}\subset \X$ forms an $m$-simplex if $\bigcap_{j=0}^m B_{r}(x_{i_j})\neq \emptyset$, where $B_r(x):=\{ y\in \T^d: \ms{dist}(x,y) \le r \}$ is the closed ball of radius $r>0$ centered at $x$.  
\end{itemize}
\end{definition}

Intrinsically,  the \v{C}ech complex possesses inclusion property $\C(\X,r) \subset \C(\X,r')$ for all $r\le r'$ and thus induces a \emph{\v{C}ech filtration} $\big( \C(\X,r) \big)_{r\ge0}$.  To analyze the homology of a \v{C}ech filtration, the authors in \cite{bobrowski:adler:2014} employed an approach based on an extension of classical Morse theory (see \cite{gershkovich:rubinstein:1997}) to the  min-type distance function $d_\mathcal X$  in \eqref{e:min.type.dist.func}. Then, for each $k\in \{1,\dots, d\}$, the change in the $k$th homology group of a \v{C}ech filtration can be characterized by the \emph{Morse critical point} with index $k$ of  $d_\X$. 
\begin{definition}  \label{def:crit.point} 
A point $c\in \R^d$ is said to be a (Morse) critical point of index $k$ if there exists a subset $\Y\subset \X$ of $k+1$ points such that 
\begin{itemize}
\item the points in $\Y$ are in general position, i.e., $\Y$ spans a $k$-dimensional simplex in $\R^d$, so that there is a unique $(k-1)$-dimensional sphere containing $\Y$. 
\item $d_\X(c) = \ms{dist}(c,y)$ for any $y\in \Y$ and $d_\X(c)<\min_{z\in \X\setminus \Y} \ms{dist}(c,z)$. 
\item $c\in \sigma (\X)$, where $\sigma (\X)$ denotes an open geometric $k$-simplex in $\R^d$ spanned by $\X$. 
\end{itemize}
\end{definition}
Whenever such $\Y$ exists, we say that $\Y$ forms a \emph{critical $k$-face} (or simplex) for which the critical point  is given by $c=c(\Y)$. Moreover,  denote by $\rho(\Y)$ its critical value, i.e., the radius of a ball spanned by $\Y$ centered at $c(\Y)$. The $0$-dimensional critical points are  $\X$ itself and hence are not of interest to us.

Recall that $\Pn$ is a homogeneous Poisson point process with intensity $n$ in the $d$-dimensional flat torus $\T^d$. Let $(a_n)_{n\ge1}$ be a sequence satisfying \eqref{e:def.an.M0}, and $(r_n(u))_{n\ge1}$ be defined as in \eqref{e:def.rn.u}. As mentioned in the Introduction,  an extra caution is needed when  the \v{C}ech filtration is defined on the torus $\T^d$. For example, if the radius $r$ is large enough,  the intersection of balls on the torus is not contractible; hence, the \v{C}ech complex $\C(\Pn,r)$ will not be homotopy equivalent to the union of balls $\bigcup_{p\in \Pn}B_{r}(p)$ (see the Nerve Lemma  in Theorem 10.7 of \cite{bjorner:1995}). 
Moreover, when $r$ is large, the notion of a critical point, as well as its critical value, is not always well-defined on the torus. To overcome this issue, we follow the convention of the previous studies  in \cite{bobrowski:2022, bobrowski:oliveira:2019},  and focus only on a ``bounded" filtration $(\C(\Pn,r))_{r\in [0,R_n]}$, where $R_n$ satisfies 
\begin{equation}  \label{e:Rn.cond}
R_n\to0,  \ \ \ \frac{r_n(u)}{R_n}\to0,  \ \ \ n\to\infty, 
\end{equation}
for all $u\in \R$. Note that \eqref{e:Rn.cond} implies $a_n/(nR_n^d)\to0$ as $n\to\infty$.

Let $k\in \{ 1,\dots,d-1 \}$ be a positive integer. For a $(k+1)$-point subset $\Y\subset \Pn$, which is in general position, we define 
\begin{equation}  \label{e:def.s.kn}
s_{k,n}(\Y,\Pn):= \one \{ \Y \text{ forms a critical } k\text{-face} \}\times \one \{ \rho(\Y)\le R_n \}, 
\end{equation}
where $\rho(\Y)$ is the critical value defined  after Definition \ref{def:crit.point}. 
Next, for $u\in \R$, we define 
\begin{align}
\begin{split}  \label{e:def.g.kn}
g_{k,n}(\Y, \Pn; u) &:= s_{k,n}(\Y, \Pn) \one \{ \rho(\Y)\ge r_n(u)  \}\\
&= \one \{ \Y \text{ forms a critical } k\text{-face} \}\times \one \big\{ \rho(\Y)\in [r_n(u), R_n] \big\}. 
\end{split}
\end{align}
Using \eqref{e:def.s.kn}, the point process of our interest can be formally defined as 
\begin{equation}  \label{e:def.eta.kn}
\eta_{k,n}:= \sum_{\Y\subset \Pn, \, |\Y|=k+1} s_{k,n}(\Y, \Pn)\,\delta_{(c(\Y), \,n\omega_d \rho(\Y)^d-a_n)}. 
\end{equation}
Notice that $\omega_d \rho(\Y)^d$ represents the volume of an open ball in $\R^d$ with radius $\rho(\Y)$ centered at $c(\Y)$. The process \eqref{e:def.eta.kn} is viewed as a random element in the space $M_p(\bbY)$ of point measures on $\bbY:=\T^d\times (-\infty,\infty]$. Defining 
$$
b_{k,n} := na_n^{k-1}e^{-a_n}, \ \ \ n\ge 1, 
$$
we consider the sequence of probability measures 
\begin{equation}  \label{e:eta.kn.prob.measure}
\big( b_{k,n}^{-1}\P(\eta_{k,n}\in \cdot), \, n\ge1 \big). 
\end{equation}
As described in the Introduction, the convergence of \eqref{e:eta.kn.prob.measure} has to be treated under   the $\M_0$-topology. 

The formal definition of $\M_0$-topology is given as follows. First, let $B_{\emptyset, r}$ denote an open ball of radius $r>0$ centered at the null measure $\emptyset$ (in terms of the vague metric). Denote  by $\M_0=\M_0(M_p(\bbY))$ the space of Borel measures on $M_p(\bbY)$, the restriction of which to $M_p(\bbY)\setminus B_{\emptyset, r}$ is finite for all $r>0$. Moreover, define $\mathcal C_0=\mathcal C_0(M_p(\bbY))$ to be the space of continuous and bounded real-valued functions on $M_p(\bbY)$ that vanish in the neighborhood of $\emptyset$. Given $\xi_n,\xi\in\M_0$, we say that $\xi_n$ converges to $\xi$ in the $\M_0$-topology, denoted as $\xi_n\to\xi$ in $\M_0$, if it holds that $\int_{\M_p(\bbY)}g(\mu)\xi_n(\dif \mu)\to\int_{\M_p(\bbY)}g(\mu)\xi(\dif \mu)$ for all $g\in\mathcal C_0$. One may refer to \cite{hult:lindskog:2006a} for more detailed discussion on $\M_0$-topology.

\begin{theorem}  \label{t:M0.eta.kn}
Let $1\le k \le d-1$. We have, as $n\to\infty$, 
\begin{equation}  \label{e:M0.eta.kn}
b_{k,n}^{-1}\P(\eta_{k,n}\in \cdot) \to \lambda_k, \ \ \text{ in } \M_0, 
\end{equation}
where 
\begin{equation}  \label{e:def.lambda.k}
\lambda_k (\cdot) = D_k \int_\bbY \one \{ \delta_{(c,u)}\in \cdot  \} \, e^{-u} \dif c \dif u, 
\end{equation}
with $D_k$ being a positive constant defined specifically at \eqref{e:def.Dk}.
\end{theorem}
Observe that $D_k^{-1}\lambda_k$ is a measure independent of $k$ and concentrated on singletons. For technical reasons, we are forced to skip the case $k = d$.
\begin{example}
We  consider the centering term
\begin{equation}  \label{e:logloglogn}
a_n = \log n + (k-1)\log \log n + \log \log \log n. 
\end{equation}
It is then easy to calculate that $b_{k,n}\sim (\log \log n)^{-1}$ as $n\to\infty$, and by Theorem \ref{t:M0.eta.kn}, 
\begin{equation}  \label{e:log.decay}
(\log \log n) \P(\eta_{k,n}\in \cdot) \to \lambda_k, \ \ \text{ in } \M_0. 
\end{equation}
Note that $(a_n)_{n\ge1}$ in \eqref{e:logloglogn} satisfies condition \eqref{e:def.an.M0}, but its growth rate is  close to the  sequence \eqref{e:def.an.Poisson.intro} of the Poisson regime (the difference  is at most of order $\mathcal O(\log \log \log n)$). As a consequence, the probability distribution of $\eta_{k,n}$ decays only logarithmically. 

In contrast, if one takes $a_n=2\log n$, then much fewer number of  critical $k$-faces will be  counted by  the process $(\eta_{k,n})_{n\ge1}$. It then follows that $b_{k,n}\sim (2\log n)^{k-1}n^{-1}$, and
$$
 (2\log n)^{-(k-1)}n\P(\eta_{k,n}\in \cdot) \to \lambda_k, \ \ \text{ in } \M_0. 
$$
In this case, the probability distribution of $\eta_{k,n}$ decays much faster than \eqref{e:log.decay}. 
\end{example}
\medskip
{In parallel, we also study the asymptotics of the sequence $(b_{k,n}^{-1}\P\circ G_{k,n}^{-1})_{n\ge 1}$, where $G_{k,n}$ is formally defined as follows. For a fixed $u_0\in \R$,
\begin{equation}  \label{e:def.G.kn}
G_{k,n}:= G_{k,n}(u_0) = \sum_{\Y\subset \Pn, \, |\Y|=k+1} g_{k,n}(\Y, \Pn; u_0), 
\end{equation}
which   counts  the number of critical points $c(\Y)\in \T^d$ with index $k$, such that $\rho(\Y)\in [r_n(u_0), R_n]$.}

Let $E:=(0,\infty]$ and $M_+(E)$ be the space of Radon measures on $E$. Define $C_K^+(E)$ to be the collection of non-negative and continuous functions on $E$ with compact support. Given $\xi_n,\xi\in M_+(E)$, we say that $\xi_n$ {converges vaguely} to $\xi$, denoted as $\xi_n \stackrel{v}{\to} \xi$ in $M_+(E)$, if it holds that $\int_{E}g(x)\xi_n(\dif x)\to\int_{E}g(x)\xi(\dif x)$ for all $g\in C_K^+(E)$. 

\begin{theorem}  \label{t:vague.eta.kn}
Under the assumptions of Theorem \ref{t:M0.eta.kn}, as $n\to\infty$, 
$$
b_{k,n}^{-1} \P(G_{k,n}\in \cdot) \stackrel{v}{\to} D_k e^{-u_0}\delta_1,  \ \ \text{ in } M_+(E), 
$$
where $u_0$ is a fixed real number as at \eqref{e:def.G.kn} and $D_k$ is a positive constant defined at \eqref{e:def.Dk}.\\
This implies that 
$$
b_{k,n}^{-1} \P(G_{k,n}\ge1) \to D_k e^{-u_0}. 
$$

\end{theorem}

From the Morse-theoretic analyses on critical faces in \cite{bobrowski:2022}, a critical $k$-face either generates a (nontrivial) $k$-dimensional cycle in the $k$th homology group of a \v{C}ech filtration or terminates a $(k-1)$-dimensional cycle of the same complex. In the former case, we call such a critical $k$-face a \emph{positive critical $k$-face}, while the latter one is called a \emph{negative critical $k$-face}. {For example, a negative critical $1$-face is nothing but an edge in the minimal spanning tree on $\Pn$ with weights being the Euclidean distance.} Now, we introduce a series of indicator functions analogous to \eqref{e:def.s.kn} and \eqref{e:def.g.kn}: for a $(k+1)$-point subset $\Y\subset \Pn$, which is in general position, 
\begin{align*}
s_{k,n}^+(\Y,\Pn) &:= \one \{  \Y \text{ forms a positive critical } k\text{-face} \}\times \one \{ \rho(\Y)\le R_n \}, \\
s_{k,n}^-(\Y,\Pn) &:= \one \{  \Y \text{ forms a negative critical } k\text{-face} \}\times \one \{ \rho(\Y)\le R_n \}, 
\end{align*}
and, for $u\in \R$, 
$$
g_{k,n}^\pm (\Y, \Pn; u) := s_{k,n}^\pm (\Y, \Pn) \one \{ \rho(\Y)\ge r_n(u) \}. 
$$
Analogously to \eqref{e:def.eta.kn} and \eqref{e:def.G.kn}, we define the point processes, 
\begin{equation}  \label{e:def.eta.kn.pm}
\eta_{k,n}^\pm:= \sum_{\Y\subset \Pn, \, |\Y|=k+1} s_{k,n}^\pm(\Y, \Pn)\,\delta_{(c(\Y), \,n\omega_d \rho(\Y)^d-a_n)} \in M_p(\bbY), 
\end{equation}
and for a fixed $u_0\in \R$, 
$$
G_{k,n}^\pm := G_{k,n}^\pm(u_0)  = \sum_{\Y\subset \Pn, \, |\Y|=k+1} g_{k,n}^\pm (\Y, \Pn; u_0). 
$$
It was shown in  Proposition 4.2 and 5.1 of  \cite{bobrowski:2022} that if $(a_n)$ satisfies condition \eqref{e:def.an.Poisson.intro}, then $(G_{k,n})_{n\ge1}$, $(G_{k,n}^+)_{n\ge1}$, and $(G_{k+1,n}^-)_{n\ge1}$ exhibit  the same phase transition  in terms of  the $k$th homological connectivity. From this point of view, the results below claim that even under the assumption \eqref{e:def.an.M0}, $(G_{k,n}^+)_{n\ge1}$ and $(G_{k+1,n}^-)_{n\ge1}$ still satisfy the same limit theorem as Theorem \ref{t:vague.eta.kn}. 
Moreover, the processes $(\eta_{k,n}^+)_{n\ge1}$ and $(\eta_{k+1,n}^-)_{n\ge1}$ also satisfy the same limit theorem as Theorem \ref{t:M0.eta.kn}. 
For the latter process, however, we can deduce the corresponding limit result only for the process with a restricted state space, i.e., 
\begin{equation}  \label{e:def.restricted.eta}
\eta_{k+1,n}^{(-,\ms{r})}:= \sum_{\Y\subset \Pn, \, |\Y|=k+2} s_{k+1,n}^-(\Y, \Pn)\,\delta_{n\omega_d \rho(\Y)^d-a_n} \in M_p\big((-\infty, \infty]\big), 
\end{equation}
for which  the locational coordinate $c(\Y)$ has been removed  from \eqref{e:def.eta.kn.pm}. 
We conjecture that  the same limit theorem  holds even for the process $(\eta_{k+1,n}^-)_{n\ge1}$ as well. However,  due to a technical complication occurring in the application of a series of results in \cite{bobrowski:2022}, we have decided not to pursue this direction in the present work. 
\begin{corollary}  \label{cor:pos.neg.faces}
Recall $\lambda_k$ as defined in \eqref{e:def.lambda.k}. \\
$(i)$ For $1\le k \le d-1$, we have 
\begin{equation}  \label{e:M0.eta.kn1}
b_{k,n}^{-1} \P(\eta_{k,n}^+\in \cdot) \to \lambda_k, \ \ \text{ in } \M_0, \ \ \ n\to\infty.
\end{equation}
$(ii)$ For $1 \le k \le d-2$, we have 
\begin{equation}  
\label{e:M0.eta.kn2}
b_{k,n}^{-1} \P(\eta_{k+1,n}^{(-,\ms{r})}\in \cdot) \to \lambda_k^{(\ms{r})}, \ \ \text{ in } \M_0, \ \ \ n\to\infty, 
\end{equation}
where $\lambda_k^{(\ms{r})}(\cdot) =  D_k \int_\bbr \one \{ \delta_{u}\in \cdot  \} \, e^{-u} \dif u, $ is the restricted version of $\lambda_k$.\\

\noindent $(iii)$ For $1\le k \le d-1$, we have  
\begin{align}
&b_{k,n}^{-1} \P(G_{k,n}^+\in \cdot) \stackrel{v}{\to}D_ke^{-u_0}\delta_1, \ \ \text{ in } M_+(E), \ \ \ n\to\infty.\label{e:vague.Gkn1}
\end{align}
$(iv)$ For $1\le k \le d-2$, we have 
\begin{equation}  \label{e:vague.Gkn2}
b_{k,n}^{-1} \P(G_{k+1,n}^-\in \cdot) \stackrel{v}{\to}D_ke^{-u_0}\delta_1, \ \ \text{ in } M_+(E), \ \ \ n\to\infty. 
\end{equation}
In particular, \eqref{e:vague.Gkn1} and \eqref{e:vague.Gkn2} imply that as $n\to\infty$, 
\begin{align*}
&b_{k,n}^{-1} \P(G_{k,n}^+\ge1) \to D_k e^{-u_0},\\
&b_{k,n}^{-1} \P(G_{k+1,n}^-\ge1) \to D_k e^{-u_0}. 
\end{align*}
\end{corollary}

The proof of the limit theorem for $(\eta_{k,n}^+)_{n\ge1}$ is based on establishing asymptotic equivalence with $(\eta_{k,n})_{n\ge1}$ by using estimates in  \cite{bobrowski:2022}. For the case of $k = 1$, we need estimates from \cite{penrose:2003}. From this we can deduce convergence of $(G_{k,n}^+)_{n\ge1}$ via the continuous mapping theorem. We then compare $(G_{k,n}^+)_{n\ge1}$ with $(G_{k+1,n}^-)_{n\ge1}$ again using estimates in \cite{bobrowski:2022}. Finally, by using the results for $(G_{k+1,n}^-)_{n\ge1}$ at many choices of $u_0$'s, we can prove the limit theorem for $(\eta_{k+1,n}^{(-,\ms{r})})_{n\ge1}$.
\section{Proofs}  \label{sec:proofs}

For the proof of Theorem \ref{t:M0.eta.kn}, we need to recall that  the \emph{Kantorovich-Rubinstein distance} between the laws of point processes $\eta_1$ and $\eta_2$ on $\bbY$ is defined as 
$$
d_{\ms{KR}} \big( \mathcal L(\eta_1),  \mathcal L(\eta_2) \big) := \sup_{h} \big| \E[h(\eta_1)] - \E[h(\eta_2)]  \big|, 
$$
where $\mathcal L(\eta_i)$ denotes a probability law of $\eta_i$, and the supremum is taken over all $1$-Lipschitz functions $h:M_p(\bbY)\to \R$, with respect to the total variation distance on $M_p(\bbY)$. Recall also that the \emph{total variation distance} between two measures $\mu_1$ and $\mu_2$ on $\bbY$ is defined as 
$$
d_{\ms{TV}}(\mu_1,\mu_2) := \sup_{A\subset \bbY} \big| \mu_1(A)-\mu_2(A) \big|. 
$$

In the previous section,  we defined the notion of a critical point $c(\Y)$, as well as its critical value $\rho(\Y)$, only when $\Y$ forms a critical face; see Definition \ref{def:crit.point}. In this section, however, one needs to extend these concepts, even when $\Y$ may not form a critical face. Specifically, let $\Y=\{ y_1,\dots,y_{k+1} \}$ be a subset of $\Pn$, which is in general position ($\Y$ does not necessarily form a critical $k$-face). Let 
$$
E(\Y):=\big\{  z\in \T^d: \| z-y_1\| = \cdots = \|z-y_{k+1} \|  \big\}
$$
be the collection of equidistant points from $\Y$. Since $\Y$ is in general position, $E(\Y)$ forms a $(d-k)$-dimensional affine plane, so that  there exists a unique point $c(\Y)\in E(\Y)$ such that $\|c(\Y)-y_1\| = \inf_{z\in E(\Y)}\|z-y_1\|$. Moreover, we define  $\rho(\Y):=\|c(\Y)-y_1\|$. 

Before commencing the proof of Theorem \ref{t:M0.eta.kn}, we define a finite and positive constant
\begin{equation}  \label{e:def.Dk}
D_k := \frac{(k!)^{d-k+1}}{(k+1)! d \omega_d^k}\, \binom{d}{k} \frac{\Omega_d}{\Omega_k \Omega_{d-k}} \int_{(S^{k-1})^{k+1}} h_k (\btheta) V_{\ms{simp}} (\btheta)^{d-k+1}\dif \btheta, 
\end{equation}
where $\Omega_j=\prod_{i=1}^j \omega_i$ ($\omega_i$ is volume of the $i$-dimensional unit ball), $S^{k-1}$ is the $(k-1)$-dimensional unit sphere, and $h_k(\btheta):=\one \{ c(\btheta)\in \sigma (\btheta) \}$ (recall that $\sigma (\btheta)$ is an open geometric $k$-simplex in $\R^d$ spanned by $\btheta$).
Further, $V_{\ms{simp}}(\btheta)$ represents the volume of $\sigma (\btheta)$. 
%a $k$-simplex spanned by $\btheta\in (S^{k-1})^{k+1}$. 

Throughout this section, denote by $C^*$ a generic positive constant, which is independent of $n$ but may vary between and within the lines. 

\subsection{Proof of Theorem \ref{t:M0.eta.kn}}

\begin{proof}
The proof exploits the ideas in Theorem 4.1 of \cite{bobrowski:schulte:yogeshwaran:2022}. Recall that $C_K^+(\bbY)$ is the collection of non-negative and continuous functions on $\bbY$ with compact support. Given $H_1, H_2\in C_K^+(\bbY)$ and $\vep_1, \vep_2>0$, define $F_{H_1,H_2,\vep_1,\vep_2}:M_p(\bbY)\to [0,1]$ by 
\begin{equation}   \label{e:def.F}
F_{H_1, H_2, \vep_1, \vep_2}(\xi)=(1-e^{-(\xi(H_1)-\epsilon_1)_{+}})(1-e^{-(\xi(H_2)-\epsilon_2)_{+}}),
\end{equation}
where $(a)_+=a$ if $a\ge 0$ and 0 otherwise, and $\xi(H_\ell)=\int_\bbY H_\ell(x)\xi(\dif x)$ for $\ell=1,2$. 
Denote $\lambda_{k,n}(\cdot):= b_{k,n}^{-1}\P(\eta_{k,n}\in \cdot)$ and recall the definition of $\lambda_k$ at \eqref{e:def.lambda.k}. According to Theorem A.2 in \cite{hult:samorodnitsky:2010}, \eqref{e:M0.eta.kn} immediately follows if one can show that 
$$
\lambda_{k,n}(F_{H_1, H_2, \vep_1, \vep_2}) \to \lambda_k (F_{H_1, H_2, \vep_1, \vep_2}),  \  \ \text{ as } n\to\infty, 
$$
for all $H_1,H_2 \in C_K^+(\bbY)$ and $\vep_1, \vep_2 >0$. 

Fix $H_1, H_2$ and $\vep_1,\vep_2$ henceforth and write $F=F_{H_1, H_2, \vep_1, \vep_2}$. Let $\zeta_{k,n}$ be a Poisson point process on $\bbY$ with intensity measure $b_{k,n} D_k e^{-u} \dif c \dif u$ for $c\in \T^d$, $u \in \R$. Then, writing $\lambda_{k,n}(F)=b_{k,n}^{-1}\E\big[  F(\eta_{k,n})\big]$, it suffices to demonstrate that as $n\to\infty$, 
\begin{equation}  \label{e:cond1.M0}
b_{k,n}^{-1} \big| \, \E\big[ F(\eta_{k,n}) \big] - \E\big[ F(\zeta_{k,n}) \big] \, \big| \to 0, 
\end{equation}
\begin{equation}  \label{e:cond2.M0}
b_{k,n}^{-1} \E\big[ F(\zeta_{k,n}) \big] \to \lambda_k(F). 
\end{equation}
We begin with proving \eqref{e:cond1.M0}. 
Since $H_\ell$ has compact support on $\bbY$, one can find $u_0\in \R$ such that 
$$
\text{supp}(H_1)\bigcup  \text{supp}(H_2) \subset \bbY_0:= \T^d \times [u_0,\infty], 
$$
where $\text{supp}(H_\ell)$ denotes the support of $H_\ell$. Thus, we may assume without loss of generality that $\eta_{k,n}$ is a point process with state space restricted to $\bbY_0$, i.e., 
\begin{equation}  \label{e:restricted.eta.kn}
\eta_{k,n} = \sum_{\Y\subset \Pn, \, |\Y|=k+1} g_{k,n}(\Y, \Pn; u_0) \delta_{(c(\Y), \, n\omega_d \rho(\Y)^d -a_n)} \in M_p(\bbY_0). 
\end{equation}
Accordingly, $\zeta_{k,n}$ can  be viewed as a Poisson point process with mean measure 
$$
(\ms{Leb}\otimes \tau_{k,n})(\dif c, \dif u):= b_{k,n} D_k e^{-u} \one \{ u\ge u_0 \} \dif c \dif u, \ \ c\in \T^d, \, u \in \R. 
$$
Since $0\le F\le 1$, it is elementary to see that 
$$
\big| F(\mu_1)-F(\mu_2)  \big| \le 2d_{\ms{TV}} (\mu_1, \mu_2), \ \ \ \mu_1, \mu_2 \in M_p(\bbY_0). 
$$
This implies that $F$ is {$2$}-Lipschitz, and hence, 
$$
\big| \, \E\big[ F(\eta_{k,n}) \big] - \E\big[ F(\zeta_{k,n}) \big] \, \big|\le 2 d_\ms{KR} \big( \mathcal L(\eta_{k,n}), \mathcal L(\zeta_{k,n})\big). 
$$
Now, \eqref{e:cond1.M0} will follow if we can prove that 
\begin{equation}  \label{e:conv.KR.dist}
b_{k,n}^{-1} d_{\ms{KR}} \big(\mathcal L (\eta_{k,n}), \mathcal L(\zeta_{k,n})\big)\to 0, \ \ \text{as } n\to\infty. 
\end{equation}
By virtue of  Theorem 4.1 in \cite{bobrowski:schulte:yogeshwaran:2022} (see also equ.~(6.5) therein), it suffices to  show that as $n\to\infty$, 
\begin{equation}  \label{e:cond1.thm4.1}
b_{k,n}^{-1} d_\ms{TV} \big( \E[ \eta_{k,n}(\cdot) ], \, \ms{Leb}\otimes \tau_{k,n} \big) \to 0, 
\end{equation}
\begin{equation}  \label{e:cond2.thm4.1}
b_{k,n}^{-1}\big\{ \text{Var}(\eta_{k,n}(\bbY_0)) -\E[\eta_{k,n}(\bbY_0)]\big\} \to0,
\end{equation}
and 
\begin{align}  
\begin{split}  \label{e:cond3.thm4.1}
&\frac{b_{k,n}^{-1}n^{2(k+1)}}{\big( (k+1)! \big)^2}\, \int_{(\T^d)^{k+1}}\int_{(\T^d)^{k+1}} \one \big\{ \B(\bx)\cap \B(\bz)\neq \emptyset \big\} \\
&\qquad \qquad \qquad\qquad \qquad \times \E\big[ g_{k,n}(\bx,\Pn+\delta_\bx; u_0) \big] \E\big[ g_{k,n}(\bz,\Pn+\delta_\bz; u_0) \big]\dif \bz \dif \bx \to 0, 
\end{split}
\end{align}
where $\B(\bx)$ denotes an open ball in $\R^d$ with radius $\rho(\bx)$ centered at $c(\bx)$. {We now prove these equations in that order.}
\medskip

\noindent \textit{Proof of \eqref{e:cond1.thm4.1}}: 
According to Lemma 2.4 in \cite{bobrowski:weinberger:2017}, 
\begin{equation}  \label{e:Lemma2.4}
\one \{ \Y \text{ forms a critical } k\text{-face} \} = h_k(\Y)  \one \big\{ \B(\Y)\cap \Pn =\emptyset\big\},
\end{equation}
where $h_k(\Y)=\one \{ c(\Y)\in \sigma (\Y) \}$.
Substituting \eqref{e:Lemma2.4} into \eqref{e:def.g.kn} with $u=u_0$, one can write 
\begin{align}
\begin{split}  \label{e:gkn.representation}
g_{k,n}(\Y, \Pn; u_0) 
&= h_k(\Y)\one \big\{ \B(\Y)\cap \Pn =\emptyset \big\} \times \one \big\{ \rho(\Y)\in [r_n(u_0), R_n] \big\}. 
\end{split}
\end{align}

Appealing to the multivariate Mecke formula for Poisson point processes (see, e.g., Chapter 4 in \cite{last:penrose:2017}) and using \eqref{e:gkn.representation}, we have, for every $A\subset \bbY_0$, 
\begin{align}
\begin{split}  \label{e:mecke1}
\E[\eta_{k,n}(A)] &= \frac{n^{k+1}}{(k+1)!}\,  \int_{(\T^d)^{k+1}} \one \big\{ (c(\bx), n\omega_d \rho(\bx)^d-a_n)\in A \big\}\E\big[ g_{k,n}(\bx,\Pn+\delta_\bx; u_0) \big] \dif \bx \\
&=\frac{n^{k+1}}{(k+1)!}\,  \int_{(\T^d)^{k+1}} \one \big\{ (c(\bx), n\omega_d \rho(\bx)^d-a_n)\in A \big\} 
h_k(\bx)\one \{ \rho(\bx) \le R_n \} e^{-n\omega_d \rho(\bx)^d}\dif \bx. 
\end{split}
\end{align}
By {a change of variable based on} the Blaschke-Petkantschin-type formula provided in Lemma C.1 of \cite{bobrowski:2022}, 
\begin{align}
\begin{split}  \label{e:mecke2}
\E[\eta_{k,n}(A)] &= \frac{D_{bp}}{(k+1)!}\, \int_{(S^{k-1})^{k+1}} h_k (\btheta) V_\ms{simp} (\btheta)^{d-k+1}\dif \btheta  \\
&\quad \times n^{k+1} \int_{\T^d}\int_0^{R_n}   \one \big\{ (c, n\omega_d \rho^d-a_n)\in A \big\} \rho^{dk-1}e^{-n\omega_d \rho^d} \dif \rho \dif c,
\end{split} 
\end{align}
where 
$$
D_{bp} =(k!)^{d-k+1} \binom{d}{k} \frac{\Omega_d}{\Omega_k \Omega_{d-k}}. 
$$
Performing the change of variable by $u=n\omega_d \rho^d-a_n$, 
\begin{align*}
&n^{k+1} \int_{\T^d}\int_0^{R_n}   \one \big\{ (c, n\omega_d \rho^d-a_n)\in A \big\} \rho^{dk-1}e^{-n\omega_d \rho^d} \dif \rho \dif c \\
&=\frac{b_{k,n}}{d\omega_d^k}\, \int_A \one \big\{ u \in (-a_n, n\omega_d R_n^d -a_n) \big\} \Big( 1+\frac{u}{a_n} \Big)^{k-1} e^{-u} \dif c \dif u. 
\end{align*}
It thus follows that 
\begin{equation} \label{e:E[eta.kn(A)]}
\E[\eta_{k,n}(A)]=D_k b_{k,n} \int_A \one \big\{ u \in (-a_n, n\omega_d R_n^d -a_n) \big\} \Big( 1+\frac{u}{a_n} \Big)^{k-1} e^{-u} \dif c \dif u, 
\end{equation}
where $D_k$ is defined in \eqref{e:def.Dk}.

Since 
$$
(\ms{Leb}\otimes \tau_{k,n})(A) = b_{k,n}D_k \int_Ae^{-u} \dif c \dif u, 
$$
we have, as $n\to\infty$, 
\begin{align*}
b_{k,n}^{-1} d_\ms{TV} \big( \E[ \eta_{k,n}(\cdot) ], \, \ms{Leb}\otimes \tau_{k,n} \big) &=b_{k,n}^{-1} \sup_{A\subset \bbY_0} \big|\, \E[\eta_{k,n}(A)]- (\ms{Leb}\otimes \tau_{k,n})(A)  \big|  \\
&\le D_k\int_{u_0}^\infty \Big|  \Big( 1+\frac{u}{a_n} \Big)^{k-1}-1 \Big| e^{-u} \dif u \\
&\quad + D_k \int_{u_0}^\infty  \one \big\{ u \notin (-a_n, n\omega_d R_n^d -a_n)\big\} e^{-u} \dif u \to 0. 
\end{align*}
The last convergence follows from the dominated convergence theorem and the assumption that $a_n/(nR_n^d)\to0$, $n\to\infty$. 
\medskip

\noindent \textit{Proof of \eqref{e:cond2.thm4.1}}: By a simple calculation, 
\begin{align*}
\E\big[ \eta_{k,n}(\bbY_0)^2 \big] &= \E\Big[ \Big( \sum_{\Y\subset \Pn, \, |\Y|=k+1} g_{k,n}(\Y, \Pn; u_0) \Big)^2 \Big] \\
&= \sum_{\ell=0}^{k+1}\E \Big[  \sum_{\substack{\Y\subset \Pn, \\ |\Y|=k+1}}\sum_{\substack{\Y'\subset \Pn, \\ |\Y'|=k+1}} g_{k,n}(\Y, \Pn; u_0) g_{k,n}(\Y', \Pn; u_0) \one \{ |\Y\cap \Y'|=\ell \}  \Big]\\
&=:\sum_{\ell=0}^{k+1}I_{\ell,n}. 
\end{align*}
One can compute $I_{\ell,n}$  explicitly using Palm theory for Poisson processes; see Section 8 in \cite{bobrowski:2022}. We will now use some of the bounds therein. Since $I_{k+1,n}=E[\eta_{k,n}(\bbY_0)]$, it follows that 
$$
\text{Var}(\eta_{k,n}(\bbY_0)) -\E[\eta_{k,n}(\bbY_0)] = \sum_{\ell=1}^k I_{\ell,n} + (I_{0,n}-I_{k+1,n}^2). 
$$
Now, for the proof of \eqref{e:cond2.thm4.1}, one needs to show that for every $\ell\in \{ 1,\dots,k \}$, 
\begin{equation}  \label{e:var.exp1}
b_{k,n}^{-1}I_{\ell,n}\to0,\ \ \ n\to\infty,
\end{equation}
and 
\begin{equation}  \label{e:var.exp2}
b_{k,n}^{-1} (I_{0,n}-I_{k+1,n}^2)\to0, \ \ \ n\to\infty. 
\end{equation}
{For the  upper bound of \eqref{e:var.exp1}, we use equ.~(8.15) in \cite{bobrowski:2022}, while  replacing  $\Lambda := n\omega_dr_n^d(u_0)$ therein with $a_n+u_0$;  this yields that for every $\ell\in \{ 1,\dots,k \}$, }
\begin{align}
\begin{split}  \label{e:upper.bdd.bn.I.ln}
b_{k,n}^{-1}I_{\ell,n} &\le C^* b_{k,n}^{-1} n(a_n+u_0)^{k-1} e^{-(a_n+u_0)} \Big\{  \vep_\ell^{k+2-\ell} (a_n+u_0)^{k+1-\ell} \\
&\qquad + \Big( \frac{\delta_\ell}{\vep_\ell} \Big)^{d-k}  \big(\delta_\ell (a_n+u_0)\big)^{k-\ell+1}  + (a_n+u_0)^{k+1-\ell} e^{-C_1 \delta_\ell (a_n+u_0)}  \Big\}, 
\end{split}
\end{align}
where $C_1\in (0,\infty)$ is some constant and 
\begin{equation}   \label{e:vep.l.delta.l}
\vep_\ell := (a_n+u_0)^{\frac{1}{2(k+2-\ell)}-1}; \ \ \ \delta_\ell := \frac{k+2-\ell}{C_1}\cdot \frac{\log (a_n+u_0)}{a_n+u_0}. 
\end{equation}
Substituting \eqref{e:vep.l.delta.l} into \eqref{e:upper.bdd.bn.I.ln}, we obtain that 
$$
b_{k,n}^{-1}I_{\ell,n} \le C^* \big\{ a_n^{-1/2} + C^* (\log a_n)^{d-\ell+1} a_n^{-\frac{d-k}{2(k+2-\ell)}} +a_n^{-1}\big\} \to 0, \ \ \ n\to\infty. 
$$
For the proof of \eqref{e:var.exp2}, we use equ.~(8.17) in \cite{bobrowski:2022} and replace  $\Lambda := n\omega_dr_n^d(u_0)$ therein with $a_n+u_0$, to obtain that 
\begin{align*}
b_{k,n}^{-1} (I_{0,n}-I_{k+1,n}^2) &\le C^* b_{k,n}^{-1} n (a_n+u_0)^{k-1} e^{-(a_n+u_0)} \\
&\quad \times \big\{ \vep_0^{d+2} (a_n+u_0)^{k+1} + (a_n+u_0)^{k+1} e^{-C_1 \vep_0 (a_n+u_0)} \big\}, 
\end{align*}
where $\vep_0:=(a_n+u_0)^{-(k+2)/(d+2)}$. It is then easy to show that
$$
b_{k,n}^{-1} (I_{0,n}-I_{k+1,n}^2)\le C^* \big\{  a_n^{-1} +a_n^{k+1} e^{-C_1 a_n^{(d-k)/(d+2)}}  \big\} \to 0, \ \ \text{as } n\to\infty. 
$$
\medskip

\noindent \textit{Proof of \eqref{e:cond3.thm4.1}}:  First fix $\bx \in (\T^d)^{k+1}$. For any $\bz\in (\T^d)^{k+1}$ with $\B(\bx)\cap \B(\bz)\neq \emptyset$, we have  $\|c(\bx)-c(\bz)\| \le \rho(\bx)+\rho(\bz) \le 2R_n$. Thus, 
$$
\one \big\{  \B(\bx)\cap \B(\bz)\neq \emptyset \big\} \le \one \big\{  c(\bz)\in B_{2R_n}(c(\bx)) \big\}. 
$$
Applying this inequality, while proceeding as in \eqref{e:mecke1} and \eqref{e:mecke2}, 
\begin{align*}
&\frac{b_{k,n}^{-1}n^{k+1}}{(k+1)!}\,\int_{(\T^d)^{k+1}} \one \big\{  \B(\bx)\cap \B(\bz)\neq \emptyset \big\} \, \E [g_{k,n}(\bz, \Pn+\delta_\bz; u_0)] \dif \bz \\
&\le \frac{b_{k,n}^{-1}n^{k+1}}{(k+1)!}\,\int_{(\T^d)^{k+1}}  \one \big\{  c(\bz)\in B_{2R_n}(c(\bx)) \big\} h_k(\bz) \one \big\{ \rho(\bz)\in [r_n(u_0), R_n] \big\} e^{-n\omega_d \rho(\bz)^d} \dif \bz  \\
&=\frac{b_{k,n}^{-1}n^{k+1}}{(k+1)!}\, D_{bp} \int_{(S^{k-1})^{k+1}} h_k (\btheta) V_\ms{simp}(\btheta)^{d-k+1} \dif \btheta \int_{B_{2R_n}(c(\bx))} \dif c \int_{r_n(u_0)}^{R_n} \rho^{dk-1}e^{-n\omega_d \rho^d} \dif \rho  \\
&=b_{k,n}^{-1} (2R_n)^d \omega_d \E[\eta_{k,n}(\bbY_0)], 
\end{align*}
where the last equality is due to \eqref{e:mecke2}  with $A=\bbY_0$. 
Because of \eqref{e:cond1.thm4.1}, the last term above is further bounded by 
$C^* b_{k,n}^{-1} R_n^d (\ms{Leb}\otimes \tau_{k,n})(\bbY_0) =C^* R_n^d$. 
Referring the obtained bound back to \eqref{e:cond3.thm4.1},  one can eventually bound  \eqref{e:cond3.thm4.1}  by 
$$
C^* R_n^d \,\frac{n^{k+1}}{(k+1)!}\, \int_{(\T^d)^{k+1}} \E [g_{k,n}(\bx, \Pn+\delta_\bx; u_0)] \dif \bx = C^* R_n^d \E[\eta_{k,n}(\bbY_0)] = o(R_n^d) \to 0, \ \ \ n\to\infty, 
$$
as desired. 
\vspace{5pt}

We will proceed to show  \eqref{e:cond2.M0}. Notice that one may express $\zeta_{k,n}=\sum_{i=1}^{N_{k,n}} \delta_{(C_i, U_i)}$, where $N_{k,n}$ is a Poisson random variable with mean $b_{k,n}D_ke^{-u_0}$ and $(C_i,U_i)$ are i.i.d.~random vectors on $\bbY$ with density $e^{-(u-u_0)}\one \{ u\ge u_0 \}\dif c \dif u$. Besides, $N_{k,n}$ {can be} taken to be independent of $(C_i,U_i)$. 
With this {construction} now available, we have 
\begin{align*}
b_{k,n}^{-1} \E[F(\zeta_{k,n})] &=b_{k,n}^{-1} \E \Big[ \prod_{\ell=1}^2 \Big( 1-e^{-\big( \sum_{i=1}^{N_{k,n}} H_\ell(C_i,U_i) -\vep_\ell\big)_+} \Big) \Big]  \\
&= b_{k,n}^{-1} \E \Big[ \prod_{\ell=1}^2 \Big( 1-e^{-(H_\ell(C_1,U_1) -\vep_\ell)_+} \Big)\one \{ N_{k,n}=1 \} \Big]  \\
&\qquad \qquad + b_{k,n}^{-1} \E \Big[ \prod_{\ell=1}^2 \Big( 1-e^{-\big( \sum_{i=1}^{N_{k,n}} H_\ell(C_i,U_i) -\vep_\ell\big)_+} \Big)\one \{ N_{k,n}\ge 2 \} \Big]  \\
&=:A_n+B_n. 
\end{align*}
By an elementary calculation, 
$$
B_n\le \bkn^{-1} \P(N_{k,n}\ge2) \le C^*\bkn\to0, \ \ \ n\to\infty. 
$$
By the independence of $N_{k,n}$ and $(C_1,U_1)$, 
\begin{align*}
A_n &= b_{k,n}^{-1} \E \Big[ \prod_{\ell=1}^2 \Big( 1-e^{-(H_\ell(C_1,U_1) -\vep_\ell)_+} \Big)\Big] \P(N_{k,n}=1) \\
&=  D_k e^{-u_0}\cdot e^{-\bkn D_k e^{-u_0}} \int_\bbY\prod_{\ell=1}^2 \Big( 1-e^{-(H_\ell(c,u)-\vep_\ell)_+} \Big) e^{-(u-u_0)}\one \{u\ge u_0  \} \dif c \dif u \\
&\to D_k \int_\bbY \prod_{\ell=1}^2 \Big( 1-e^{-(H_\ell(c,u)-\vep_\ell)_+} \Big) e^{-u} \dif c \dif u = \lambda_k(F), 
\end{align*}
as required. 
\end{proof}
\medskip

\subsection{Proof of Theorem \ref{t:vague.eta.kn}}
\begin{proof}
By restricting the state space from $\bbY$ to $\bbY_0=\T^d \times [u_0,\infty]$, we can establish  $\M_0$-convergence  analogous to Theorem \ref{t:M0.eta.kn}. Namely,  as $n\to\infty$, 
$$%\label{e:restricted.M0.conv}
b_{k,n}^{-1} \P(\eta_{k,n}\in \cdot) \to \lambda_k, \ \ \text{ in } \M_0. 
$$
Here,  $\eta_{k,n}$ is defined as in \eqref{e:restricted.eta.kn} due to the restriction of the state space, while  the limiting measure is written as 
$$
\lambda_k(\cdot) =D_k \int_{\bbY_0} \one \{ \delta_{(c,u)}\in \cdot \} e^{-u} \dif c \dif u. 
$$

Define a continuous map $V:M_p(\bbY_0)\to \bbn:=\{ 0,1,2,\dots \}$ by $V(\xi)=\xi(\bbY_0)$, where  $M_p(\bbY_0)$ is equipped with vague topology, and  $\bbn$ is equipped with the discrete topology. It then follows from Theorem 2.5 in \cite{hult:lindskog:2006a} that 
$$
\lambda_{k,n}\circ V^{-1} \to \lambda_k \circ V^{-1}, \ \ \text{ in } \M_0(M_p(\bbn)), \ \ \ n\to\infty; 
$$
equivalently, 
\begin{equation}  \label{e:M0.to.vague}
b_{k,n}^{-1}\P(G_{k,n}\in \cdot) \to D_k e^{-u_0} \delta_1, \ \ \text{ in } \M_0(M_p(\bbn)), \ \ \ n\to\infty. 
\end{equation}
For every $x>0$, the indicator $\one_{[x,\infty)}(\cdot): \bbn\to \{ 0,1\}$ is {bounded and continuous} with respect to the discrete topology. Moreover, it vanishes in the neighborhood of $0$ (i.e., the origin of $\bbn$). In conclusion, $\one_{[x,\infty)}(\cdot) \in \mathcal C_0$, and thus, \eqref{e:M0.to.vague} implies that for every $x>0$, 
$$
\bkn^{-1} \P(G_{k,n}\ge x) \to D_k e^{-u_0}\delta_1( [x, \infty)),  \ \ \ n\to\infty.
$$
Now, Lemma 6.1 in \cite{resnick:2007} concludes the desired vague convergence in Theorem \ref{t:vague.eta.kn}. 
\end{proof}
\medskip

\subsection{Proof of Corollary \ref{cor:pos.neg.faces}}

\begin{proof}[Proof of \eqref{e:M0.eta.kn1}]
\noindent \textsc{\textit{Part I}: $2 \leq k \leq d-1$.} \\
It suffices to show that as $n\to\infty$, 
$$
b_{k,n}^{-1} \big| \, \E\big[ F(\eta_{k,n}) \big] - \E\big[ F(\eta_{k,n}^+) \big] \, \big| \to 0, 
$$
where $F$ is defined at \eqref{e:def.F} (subscripts are all omitted). Since $\bigcup_{\ell=1}^2 \text{supp}(H_\ell)\subset \bbY_0 = \T^d \times [u_0,\infty]$ for some $u_0\in \R$, we can reformulate $\eta_{k,n}^\pm$ as 
\begin{equation}  \label{e:restricted.eta.kn.pm}
\eta_{k,n}^\pm = \sum_{\Y\subset \Pn, \, |\Y|=k+1} g_{k,n}^\pm(\Y, \Pn; u_0) \delta_{(c(\Y), \, n\omega_d \rho(\Y)^d -a_n)}\in M_p(\bbY_0). 
\end{equation}
As $0\le F\le 1$, we have that 
$$
b_{k,n}^{-1} \big| \, \E\big[ F(\eta_{k,n}) \big] - \E\big[ F(\eta_{k,n}^+) \big] \, \big| \le \bkn^{-1} \P(\eta_{k,n}\neq \eta_{k,n}^+). 
$$
By definition, $\eta_{k,n}=\eta_{k,n}^++\eta_{k,n}^-$; hence, $\eta_{k,n}\neq \eta_{k,n}^+$ implies that $G_{k,n}^-=\eta_{k,n}^-(\bbY_0)\ge1$. Now, by Markov's inequality, 
$$
\bkn^{-1} \P(\eta_{k,n}\neq \eta_{k,n}^+) \le \bkn^{-1} \P(G_{k,n}^-\ge 1)\le \bkn^{-1}\E[G_{k,n}^-]. 
$$
By the proof of Proposition 5.1 in \cite{bobrowski:2022} (see the equation at page 744, line 3 therein), 
\begin{equation}  \label{e:prop5.1}
 \bkn^{-1}\E[G_{k,n}^-]\le C^* \bkn^{-1} n (a_n+u_0)^{k-2} e^{-(a_n+u_0)} \le C^* a_n^{-1}\to 0, \ \ \text{as } n\to\infty, 
\end{equation}
for which we have replaced ``$\Lambda$" in \cite{bobrowski:2022} with $a_n+u_0$. This completes the proof of Part I.
\medskip

\noindent \textsc{\textit{Part II}: $k = 1$.} \\
Our proof uses  ideas from Chapter 13 of  \cite{penrose:2003}. As in Part I, $\eta_{1,n}$ and $\eta_{1,n}^+$ can be given in the form of \eqref{e:restricted.eta.kn} and \eqref{e:restricted.eta.kn.pm} for a fixed $u_0\in \R$. For $r>0$, let $G(\Pn, r)$ be a random geometric graph on a Poisson point process $\Pn$ with edges $\{ x,y \}$ for all pairs $x,y\in \Pn$ with $\|x-y\|\le r$. Since $\eta_{1,n}^-$ corresponds to the edges of the minimal spanning tree on $\Pn$, $G_{1,n}^- \geq 1$ implies that $G(\Pn, 2r_n(u_0))$ is not connected, and thus, by arguing as in Part I we have that
$$
\P(\eta_{1,n}\neq \eta_{1,n}^+)\le \P\big(G_{1,n}^- \geq 1\big) \le \P\big( G(\Pn, 2r_n(u_0)) \text{ is not connected}\big).
$$
Hence, it  suffices to demonstrate that
\begin{equation}  \label{e:connectivity.G}
b_{1,n}^{-1} \P\big( G(\Pn, 2r_n(u_0)) \text{ is not connected}\big) \to 0, \ \ \ n\to\infty. 
\end{equation}

Before continuing, we shall introduce a few required notions. Given a graph $G$ with vertex set $V$, a non-empty subset $U\subset V$ is said to be a separating set for $G$ if none of the vertices in $V\setminus U$ are adjacent to $U$. Moreover, a pair of non-empty disjoint sets $U\subset V$, $W\subset V$ is called a separating pair for $G$, if $(i)$ the subgraph of $G$ induced by $U$ is connected and the same holds for $W$, $(ii)$ none of the elements in $U$ are adjacent to any element in $W$, and $(iii)$ none of the elements of $V\setminus (U\cup W)$ are adjacent to $U\cup W$.  Using these notions and given $K>0$, we define $E_n(K)$ as the event that there exists a separating set $U$ for $G(\Pn, 2r_n(u_0))$ with at least two elements, such that the diameter of $U$ is less than $2Kr_n(u_0)$. Further, denote by $F_n(K)$ the event that there exists a separating pair $(U,W)$ for $G(\Pn, 2r_n(u_0))$, so that the diameters of $U$ and $W$ both exceed $2Kr_n(u_0)$. 

Now, one can upper bound \eqref{e:connectivity.G} as 
\begin{equation}  \label{e:connectivity.split}
b_{1,n}^{-1} \P\big( G(\Pn, 2r_n(u_0)) \text{ is not connected}\big) \le b_{1,n}^{-1}\P(A_n) + b_{1,n}^{-1}\P(B_n), 
\end{equation}
where $A_n$ is the event that there exists at least one isolated vertex in $G(\Pn, 2r_n(u_0))$, and $B_n$ is the event for which there are no isolated vertices but $G(\Pn,2r_n(u_0))$ contains \emph{multiple} connected components (of size at least $2$). Now Markov's inequality and the Mecke formula for Poisson point processes yield that 
$$
b_{1,n}^{-1} \P(A_n) \le b_{1,n}^{-1}n \int_{\T^d} \P \big( \Pn(B_{2r_n(u_0)}(x))=0 \big)\dif x =e^{a_n-2^d(a_n+u_0)} \to 0, \ \ \text{as } n\to\infty. 
$$
In order to handle the remaining term in \eqref{e:connectivity.split}, it is sufficient to prove the following results. \\
$(i)$ There exists $K>0$ such that 
\begin{equation}  \label{e:connectivity.cond1}
b_{1,n}^{-1}\P(F_n(K))\to0, \ \ \ n\to\infty. 
\end{equation}
$(ii)$ For $K$ as above, there exists $K_1\in (0,K)$ such that 
\begin{align}  
&b_{1,n}^{-1}\P(E_n(K_1))\to0, \ \ \ n\to\infty,  \label{e:connectivity.cond2} \\
&b_{1,n}^{-1}\P\big(E_n(K)\setminus E_n(K_1)\big)\to0, \ \ \ n\to\infty. \label{e:connectivity.cond3} 
\end{align}
Indeed, if one can  show $(i)$ and $(ii)$ above, it follows from the proof of Theorem 13.10 in \cite{penrose:2003} and Lemma 13.1 therein that as $n\to\infty$, 
\begin{align*}
b_{1,n}^{-1} \P(B_n) &\le  b_{1,n}^{-1}\P\big(E_n(K)\big) + b_{1,n}^{-1}\P(F_n(K)) \\
&= b_{1,n}^{-1}\P(E_n(K_1)) + b_{1,n}^{-1}\P\big(E_n(K)\setminus E_n(K_1)\big) + b_{1,n}^{-1}\P(F_n(K)) \to 0. 
\end{align*}
\noindent \textit{Proof of \eqref{e:connectivity.cond1}}: By the proof of Proposition 13.13  in  \cite{penrose:2003}, there exist constants $\vep>0$, $\gamma>0$ such that we have, for $K>0$ and large enough $n$, 
\begin{align*}
\P(F_n(K)) &\le C^* \sum_{i\ge K/\vep} (2r_n(u_0))^{-2d} ie^{\gamma i} \cdot e^{-i \vep^d (2r_n(u_0))^d n}\le C^* \frac{n^2}{a_n^2} \sum_{i\ge K/\vep} e^{-((2\vep)^d(a_n+u_0)/\omega_d-\gamma')i}, 
\end{align*}
for some $\gamma'>0$ with $ie^{\gamma i} \le e^{\gamma' i}$ for all $i$. The last expression is further bounded by \\
$C^*n^2 a_n^{-2}e^{-2^d\vep^{d-1}Ka_n/\omega_d}$; thus, 
$$
b_{1,n}^{-1} \P(F_n(K)) \le C^* ne^{-(2^d\vep^{d-1}K/\omega_d-1)a_n}\le C^* n^{2-2^d\vep^{d-1}K/\omega_d}. 
$$
If one chooses $K>\omega_d/(2\vep)^{d-1}$, then \eqref{e:connectivity.cond1} follows as desired. 
\medskip

\noindent \textit{Proof of \eqref{e:connectivity.cond2}}: Having fixed $K$ as above, we next show \eqref{e:connectivity.cond2}. It follows from equ.~at page 297, line -3 in \cite{penrose:2003} and the proof of Lemma 13.5 in \cite{penrose:2003} that there exists $K_1\in (0,K)$ such that 
\begin{align*}
\P(E_n(K_1)) &\le C^* n (nr_n(u_0)^d)^{1-d} \int_{\T^d}\P\big( \Pn(B_{2r_n(u_0)}(x))=0 \big)\dif x \\
&= C^* n \Big(\frac{a_n+u_0}{\omega_d}  \Big)^{1-d} e^{-2^d (a_n+u_0)}  \le C^*na_n^{1-d} e^{-2^d a_n}. 
\end{align*}
Therefore, $b_{1,n}^{-1}\P(E_n(K_1)) \le C^* a_n^{1-d} e^{-(1+2^d)a_n}\to 0$ as $n\to\infty$. 
\medskip

\noindent \textit{Proof of \eqref{e:connectivity.cond3}}: By equ.~(13.43) in \cite{penrose:2003} and the last line of the proof of Lemma 13.16 therein, there exists $\xi>0$ such that  
$$
\P\big( E_n(K)\setminus E_n(K_1) \big) \le C^*n\int_{\T^d}\P\big( \Pn(B_{2r_n(u_0)}(x))=0 \big) e^{-\frac{\xi n}{2}(2r_n(u_0))^d} \dif x \le C^* ne^{-2a_n-2^{d-1}\xi a_n/\omega_d}. 
$$
This implies that as $n\to\infty$, 
$$
b_{1,n}^{-1} \P\big( E_n(K)\setminus E_n(K_1) \big) \le C^* e^{-a_n-2^{d-1}\xi a_n/\omega_d} \to 0. 
$$
\end{proof}

\begin{proof}[Proof of \eqref{e:vague.Gkn1}] 
One can exploit the same proof strategy as in Theorem \ref{t:vague.eta.kn}, with $\eta_{k,n}$ replaced by $\eta^+_{k,n}$ and \eqref{e:M0.eta.kn} replaced by \eqref{e:M0.eta.kn1}. 
%Suppose first that $2\le k \le d-1$.
%Because of Theorem \ref{t:vague.eta.kn} and  Lemma 6.1 of  \cite{resnick:2007}, it is sufficient to verify that for every $x>0$, 
%$$
%\bkn^{-1} \big\{ \P(G_{k,n}\ge x) - \P(G_{k,n}^+\ge x) \big\}\to0,  \ \ \text{as } n\to\infty. 
%$$
%Since $G_{k,n}=G_{k,n}^++G_{k,n}^-$, it follows from Markov's inequality and \eqref{e:prop5.1} that 
%\begin{align*}
%0\le \bkn^{-1} \big\{ \P(G_{k,n}\ge x) - \P(G_{k,n}^+\ge x) \big\} &\le \bkn^{-1} \P(G_{k,n}^-\ge 1) \le \bkn^{-1} \E[G_{k,n}^-] \le C^* a_n^{-1}\to 0. 
%\end{align*}

%For the remaining case of $k=1$, we  observe  that for every $x>0$, 
%$$
%0 \le \P(G_{1,n}\ge x) - \P(G_{1,n}^+\ge x) \le \P\big( G(\Pn, 2r_n(u_0)) \text{ is connected}\big). 
%$$
%The proof can now be completed by \eqref{e:connectivity.G}.
\end{proof}

\begin{proof}[Proof of \eqref{e:vague.Gkn2}]
By Lemma 6.1 in \cite{resnick:2007} and \eqref{e:vague.Gkn1}, we only need to prove that for every $x>0$, 
$$
\bkn^{-1} \big\{ \P(G_{k+1,n}^-\ge x) - \P(G_{k,n}^+\ge x)  \big\} \to 0, \ \ \ n\to\infty. 
$$
It thus suffices to show that as $n\to\infty$, 
\begin{align}
\bkn^{-1} \P(G_{k,n}^+>G_{k+1,n}^-) \to 0, \label{e:pos.neg.diff1} \\
\bkn^{-1} \P(G_{k+1,n}^- > G_{k,n}^+) \to 0. \label{e:pos.neg.diff2} 
\end{align}
We begin with proving  \eqref{e:pos.neg.diff1}.  Let $i_*:H_k\big( \bigcup_{p\in \Pn} B_{r_n(u_0)}(p) \big)\to H_k(\T^d)$ be a map induced by the inclusion $i:\bigcup_{p\in \Pn} B_{r_n(u_0)}(p) \hookrightarrow \T^d$,  where $H_k(\cdot)$ represents the $k$th homology group. By the Nerve Lemma (see, e.g., Theorem 10.7 of \cite{bjorner:1995}), there exists an isomorphism $f_*:H_k\big( \C (\Pn, r_n(u_0))\big)\to H_k\big( \bigcup_{p\in \Pn} B_{r_n(u_0)}(p) \big)$. Now, we can define the map 
$$
j_*:=i_*\circ f_*: H_k\big( \C (\Pn, r_n(u_0))\big)\to H_k(\T^d). 
$$
Suppose that $j_*$ is surjective; that is, $\text{Im}(j_*)=H_k(\T^d)$. Then, all the positive critical $k$-faces will be eventually terminated by a ``matching" negative critical $(k+1)$-face, which in turn means that $G_{k,n}^+\le G_{k+1,n}^-$. 

For   $0<r_1<r_2<\infty$  and a subset $\Y\subset \Pn$ with $|\Y|=d+1$, let $A_{r_1,r_2}(\Y)$ be the closure of an annulus $B_{r_2}(c(\Y))\setminus B_{r_1}(c(\Y))$, i.e., the closed annulus centered at the critical point $c(\Y)$. We then define 
$$
\hat G_{d,n}:= \sum_{\Y\subset \Pn, \, |\Y|=d+1} g_{d,n}(\Y, \Pn; u_0) \one \Big\{ A_{r_n(u_0), 4r_n(u_0)} (\Y)\not\subset \bigcup_{p\in \Pn} B_{r_n(u_0)}(p) \Big\}. 
$$
According to the proof of Lemma 5.7 in \cite{bobrowski:2019}, as well as Lemma 5.9 therein, one can see that if $\hat G_{d,n}=0$, then $j_*$ becomes surjective.  
Combining these observations,  it now follows that 
$$
\bkn^{-1} \P(G_{k,n}^+>G_{k+1,n}^-) \le \bkn^{-1} \P(\hat G_{d,n}\ge1) \le \bkn^{-1} \E[\hat G_{d,n}]. 
$$
Finally, appealing to the last equation in the proof of Lemma 5.8 in \cite{bobrowski:2019}, 
$$
\bkn^{-1} \E[\hat G_{d,n}]\le C^* \bkn^{-1} n (a_n+u_0)^{d-1} e^{-(a_n+u_0)(1+C^*)} \le C^*a_n^{d-k} e^{-C^*a_n}\to0, \ \ \ n\to\infty, 
$$
as desired.

Next, we focus our attention to  \eqref{e:pos.neg.diff2}. To estimate the probability in \eqref{e:pos.neg.diff2}, we need to refer to the detailed discussion on the structure of critical $k$-faces, provided in Section 7 of \cite{bobrowski:2022}. 
Below, we introduce some additional concepts and notation, which we try to keep as consistent as possible with those in \cite{bobrowski:2022}. 

Suppose in the sequel that a $(k+1)$-point subset $\Y$ of $\Pn$ is  in general position. Let $\Pi(\Y)$ be the unique linear $k$-plane centered at  $c(\Y)$ containing $\Y$, and $S^{k-1}$ be the unit sphere in $\Pi(\Y)$ (centered at $c(\Y)$). Denote by $\theta(\Y)=\{ \theta_1(\Y),\dots,\theta_{k+1}(\Y) \}$ the spherical coordinates of $\Y$ in $S^{k-1}$, and define $\hat \theta_i(\Y):=\theta(\Y)\setminus \{  \theta_i(\Y)\}$. Let 
$$
\phi_1(\Y) := \min_{1\le i \le k+1} \big\|  c(\hat \theta_i(\Y))\big\|
$$
be the   scaled distance  via  $\rho(\Y)$, between the center $c(\Y)$ and the nearest $(k-1)$-face of $\Y$. Moreover, denote by $\hat \Y_\text{min}$ such a nearest $(k-1)$-face of $\Y$ from the center $c(\Y)$. Now, let $\hat \Y_i$, $i=1,\dots,k$ be the remaining $(k-1)$-faces of $\Y$ \emph{except for} $\hat \Y_{\text{min}}$. We then define 
$$
\phi_2(\Y):= \min_{1\le i \le k} \rho(\Y)^{-1}\inf_{z\in \Pi(\hat \Y_i)}\big\|  c(\hat \Y_\text{min}) -z\big\|, 
$$
which is the distance scaled by $\rho(\Y)$, between $c(\hat \Y_{\text{min}})$ and the nearest $(k-1)$-face of $\Y$ (except for $\hat \Y_\text{min}$). 
Next, let 
$$
\hat \rho(\Y) := \rho(\Y)+ \|c(\hat \Y_\text{min}) -c(\Y)  \|, 
$$
and 
$$
\hat \B(\Y):= B_{\hat \rho(\Y)}\big(  c(\hat \Y_\text{min})\big)\setminus \B(\Y), 
$$
where $\B(\Y)$ is defined at \eqref{e:cond3.thm4.1}. Let $\hat \Pi(\Y)$ be the $(d-1)$-dimensional affine plane containing $c(\hat \Y_\text{min})$ and orthogonal to the line through $c(\Y)$ and $c(\hat \Y_\text{min})$. Given $\alpha>0$, define 
\begin{equation}  \label{e:hat.B.alpha}
\hat \B_\alpha(\Y):= \big\{  y\in \hat \B(\Y): \inf_{z\in \hat \Pi(\Y)}\| y-z \|\ge \alpha \rho(\Y) \big\}. 
\end{equation}
Since $\Pi(\hat \Y_\text{min}) \subset \hat \Pi (\Y)$, the set \eqref{e:hat.B.alpha} contains  points in $\hat \B(\Y)$ that are distant at least $\alpha \rho(\Y)$ from the plane $\Pi(\hat \Y_\text{min})$. 

Before proceeding, we define 
$$
G_{k+1,n}^{(-,+)} :=\sum_{\Y\subset \Pn, \, |\Y|=k+2} g_{k+1,n}^- (\Y, \Pn; u_0) \one \{  \hat \Y_\text{min} \text{ is a positive critical } k\text{-face} \}, 
$$
which is the number of negative critical $(k+1)$-faces whose critical values are between $r_n(u_0)$ and $R_n$, so that $\hat \Y_\text{min}$ forms a positive critical $k$-face. Because of equ.~(7.11) in \cite{bobrowski:2022}, it holds that 
\begin{equation*}  %\label{e:G--G-+}
G_{k+1,n}^- - G_{k+1,n}^{(-,+)} \le \sum_{i=1}^4 G_{k+1,n}^{(i)}, 
\end{equation*}
where $G_{k+1,n}^{(i)}$ above  are  defined  respectively as 
\begin{align*}
G_{k+1,n}^{(1)} &= \sum_{\Y\subset \Pn, \, |\Y|=k+2} g_{k+1,n}^- (\Y, \Pn; u_0)\one \{ \phi_1(\Y)>\vep_1 \}, \\
G_{k+1,n}^{(2)} &= \sum_{\Y\subset \Pn, \, |\Y|=k+2} g_{k+1,n}^- (\Y, \Pn; u_0)\one \{ \phi_1(\Y)\le \vep_1, \, \phi_2(\Y)\le \vep_2 \}, \\
G_{k+1,n}^{(3)} &= \sum_{\Y\subset \Pn, \, |\Y|=k+2} g_{k+1,n}^- (\Y, \Pn; u_0)\one \big\{ \phi_1(\Y)\le \vep_1, \, \big( \hat \B(\Y)\setminus \hat \B_{\vep_3}(\Y) \big)\cap \Pn\neq \emptyset \big\}, \\
G_{k+1,n}^{(4)} &= \sum_{\Y\subset \Pn, \, |\Y|=k+2} g_{k+1,n}^- (\Y, \Pn; u_0)\one \{ \phi_1(\Y)\le \vep_1, \, \phi_2(\Y)> \vep_2, \, \hat \B_{\vep_3}(\Y)\cap \Pn \neq \emptyset \}, \\
\end{align*}
with 
$$
\vep_1:= \frac{2}{D_{k,2}}\cdot \frac{\log (a_n+u_0)}{a_n+u_0};  \ \ \vep_2:= \big( \log (a_n+u_0) \big)^{-4}; \ \ \vep_3:=\vep_1^{2/3} 
$$
($D_{k,2}$ is a positive constant introduced in the proof of Lemma 5.6 in \cite{bobrowski:2022}). Although $\vep_1$ and $\vep_3$ above are defined analogously to equ.~(7.12) in \cite{bobrowski:2022},  the definition of $\vep_2$ is  different from equ.~(7.12) of \cite{bobrowski:2022}. 

Furthermore, define 
$$
G_{k+1,n}^{(5)} = \sum_{\Y\subset \Pn, \, |\Y|=k+2} g_{k+1,n}^- (\Y, \Pn; u_0)\one \{ \phi_1(\Y)\le \vep_1, \, \rho(\Y)>r_n(u_0), \, \rho(\hat \Y_\text{min}) < r_n(u_0) \big\}. 
$$
From the discussion in Part II of the proof of Proposition 7.1 in \cite{bobrowski:2022}, it is known  that  $G_{k+1,n}^{(i)}=0$ for $i=1,\dots,5$ implies  that $G_{k+1,n}^{(-,+)}\le G_{k,n}^+$. In conclusion, one can see that 
\begin{align*}
\bkn^{-1} \P(G_{k+1,n}^- >G_{k,n}^+) &\le \bkn^{-1} \P\Big( \sum_{i=1}^4 G_{k+1,n}^{(i)} + G_{k+1,n}^{(-,+)} >G_{k,n}^+ \Big)  \\
&\le \bkn^{-1} \P\Big( G_{k+1,n}^{(-,+)}>G_{k,n}^+, \, \sum_{i=1}^5 G_{k+1,n}^{(i)} =0 \Big) + \bkn^{-1} \P \Big( \sum_{i=1}^5 G_{k+1,n}^{(i)}\ge 1 \Big) \\
&\le \sum_{i=1}^5 \bkn^{-1} \E[G_{k+1,n}^{(i)}]. 
\end{align*}

From this analysis, it is  enough to show that for every $i\in\{1,\dots,5\}$, 
\begin{equation}  \label{e:conv.G.(i)}
\bkn^{-1} \E[G_{k+1,n}^{(i)}]\to0, \ \ \text{ as } n\to\infty. 
\end{equation}
The proof of \eqref{e:conv.G.(i)} is highly related to  Lemmas 7.9--7.13 in \cite{bobrowski:2022}. 
First, by replacing ``$\Lambda$" in \cite{bobrowski:2022} with $a_n+u_0$, the proof of Lemma 7.9 in \cite{bobrowski:2022} ensures  that 
$$
\bkn^{-1} \E[G_{k+1,n}^{(1)}] \le C^* \bkn^{-1} n(a_n+u_0)^{k-2} e^{-(a_n+u_0)} \le C^* a_n^{-1}\to0, \ \ \ n\to\infty. 
$$
Subsequently,  proceeding as in the  proof of Lemma 7.10 in \cite{bobrowski:2022}, 
\begin{align*}
\bkn^{-1} \E[G_{k+1,n}^{(2)}] &\le C^* \bkn^{-1} n (a_n+u_0)^k e^{-(a_n+u_0)} \Big\{ \frac{\vep_1\vep_{12}}{\sqrt{1-\vep_1^2}} + \frac{\vep_1(\vep_1+\vep_2)}{\vep_{12}}\Big\}  \\
&\le C^*a_n \Big\{ \frac{\vep_1\vep_{12}}{\sqrt{1-\vep_1^2}} + \frac{\vep_1(\vep_1+\vep_2)}{\vep_{12}}\Big\}, 
\end{align*}
where $\vep_{12}:=(\log (a_n+u_0))^{-2}$. 
Since $\vep_1\to0$ and $\vep_1/\vep_2\to 0$ as $n\to\infty$, one can see that $\bkn^{-1} \E[G_{k+1,n}^{(2)}]$ $\le  C^* (\log a_n)^{-1}\to0$ as $n\to\infty$. Next, it follows from the proof of Lemma 7.11 in \cite{bobrowski:2022} that 
\begin{align*}
\bkn^{-1} \E[G_{k+1,n}^{(3)}] \le C^* \bkn^{-1} \vep_3^4 n (a_n+u_0)^{k+1} e^{-(a_n+u_0)}\le C^* (\log a_n)^{8/3}a_n^{-2/3}\to 0, \ \ \ n\to\infty. 
\end{align*}
Furthermore, by the proof of Lemma 7.12 in \cite{bobrowski:2022}, 
$$
\bkn^{-1} \E[G_{k+1,n}^{(4)}] \le C^* \bkn^{-1} \vep_0^{-d} \vep_1 n (a_n+u_0)^k e^{-(a_n+u_0) (1+C^* \vep_0)}, 
$$
where $\vep_0:=(a_n+u_0)^{-3/4}$. Thus, we have 
$$
\bkn^{-1} \E[G_{k+1,n}^{(4)}] \le C^* a_n^{3d/4} (\log a_n) e^{-C^*a_n^{1/4}}\to 0, \ \ \ n\to\infty. 
$$
Finally, the proof of Lemma 7.13 in \cite{bobrowski:2022} concludes that as $n\to\infty$, 
\begin{align*}
\bkn^{-1} \E[G_{k+1,n}^{(5)}] &\le C^* \bkn^{-1} \vep_1 n (a_n+u_0)^k e^{-(a_n+u_0)} \big\{ d\vep_1^2 (a_n+u_0) + o(\vep_1^2 (a_n+u_0)) \big\} \\
&\le C^*\bkn^{-1} \vep_1^3 n (a_n+u_0)^{k+1} e^{-(a_n+u_0)} \le C^*(\log a_n)^3 a_n^{-1} \to0. 
\end{align*}
Now, \eqref{e:conv.G.(i)} has been established.  
\end{proof}

\begin{proof}[Proof of \eqref{e:M0.eta.kn2}]
Analogously to \eqref{e:def.restricted.eta}, we define the process $(\eta_{k,n}^{(+,\ms{r})})_{n\ge1}$ by dropping the locational coordinate $c(\Y)$ from \eqref{e:def.eta.kn.pm}. Note that by continuous mapping theorem, the restricted process $(\eta_{k,n}^{(+,\ms{r})})_{n\ge1}$ satisfies $\M_0$ convergence analogous to \eqref{e:M0.eta.kn1} with the limit measure $\lambda_k^{(r)}$ defined at \eqref{e:M0.eta.kn2}. 
Define $F = F_{H_1,H_2,\epsilon_1,\epsilon_2}$ as in \eqref{e:def.F}, for which the domain of $F$ is restricted to $M_p\big( (-\infty,\infty] \big)$.  Then, to complete the proof, we need to show that
\begin{equation}
\label{e:F.etakpm1}
b_{k,n}^{-1}\big|\E[F(\eta_{k,n}^{(+,\ms{r})})] - \E[F(\eta_{k+1,n}^{(-,\ms{r})})]\big| \to 0.
\end{equation}
Again the compact support of $H_{\ell}$ on $(-\infty,\infty]$ gives us $u_0\in \R$ such that 
$$
\text{supp}(H_1)\bigcup  \text{supp}(H_2) \subset  [u_0,\infty]. 
$$
From \eqref{e:vague.Gkn1} and \eqref{e:vague.Gkn2}, we have  $b_{k,n}^{-1}\P(G_{k,n}^+(u_0) \geq 2) \to 0$ and  $b_{k,n}^{-1}\P(G_{k+1,n}^-(u_0) \geq 2) \to 0$ as $n \to \infty$. Because of $0\le F\le 1$, \eqref{e:pos.neg.diff1}, \eqref{e:pos.neg.diff2}, one can obtain \eqref{e:F.etakpm1}  if we show that 
\begin{equation}
\label{e:F.etakpm2}
b_{k,n}^{-1}\E\big[|F(\eta_{k,n}^{(+,\ms{r})}) - F(\eta_{k+1,n}^{(-,\ms{r})})|\one\{  G_{k,n}^+(u_0) =G_{k+1,n}^-(u_0)= 1 \}\big] \to 0, \ \ \ n\to\infty. 
\end{equation}
Fix $\delta > 0$. Again from \eqref{e:vague.Gkn1} and \eqref{e:vague.Gkn2}, we can choose $u_1 > u_0$ so large  that 
$$ \lim_{n \to \infty}b_{k,n}^{-1}\P(G_{k,n}^+(u_1) \geq 1)\le \delta, \ \ \text{ and } \ \ 
\lim_{n \to \infty}b_{k,n}^{-1}\P(G_{k+1,n}^-(u_1) \geq 1) \leq \delta.$$
Define $A_n := \{  G_{k,n}^+(u_0) =G_{k+1,n}^-(u_0)= 1,  \, G_{k,n}^+(u_1) =G_{k+1,n}^-(u_1)= 0  \}$; then, \eqref{e:F.etakpm2} follows if we can show that
\begin{equation}
\label{e:Fdiff_1An}   
\lim_{n \to \infty} b_{k,n}^{-1}\E\big[|F(\eta_{k,n}^{(+,\ms{r})}) - F(\eta_{k+1,n}^{(-,\ms{r})})|\one_{A_n}\big] = 0.
\end{equation}
Observing that 
$$
\big|F(\eta_{k,n}^{(+,\ms{r})}) - F(\eta_{k+1,n}^{(-,\ms{r})})\big| \le 2 \sum_{\ell=1}^2 |\eta_{k,n}^{(+,\ms{r})}(H_\ell) - \eta_{k+1,n}^{(-,\ms{r})}(H_\ell)|, 
$$
\eqref{e:Fdiff_1An} is implied by 
$$
\lim_{n\to\infty} b_{k,n}^{-1} \E\big[|\eta_{k,n}^{(+,\ms{r})}(H_\ell) - \eta_{k+1,n}^{(-,\ms{r})}(H_\ell)|\one_{A_n}\big] = 0, 
$$
for each $\ell=1,2$. 

Fix $\ell\in\{ 1,2 \}$. If $A_n$ holds, we may assume, without loss of generality, that $H_\ell$ is supported on $[u_0, u_1]$. Moreover, under $A_n$, there exist random variables $X_n, X_n'\in [u_0,u_1]$ such that $\eta_{k,n}^{(+,\ms{r})}(H_\ell)=H_\ell(X_n)$ and $\eta_{k+1,n}^{(-,\ms{r})}(H_\ell)=H_\ell(X_n')$. Since $H_\ell$ is uniformly continuous on $[u_0,u_1]$, for every $\delta > 0$, there exists $\delta_0>0$ such that $\big|  H_\ell(X_n)-H_\ell(X_n')\big|\le \delta$ whenever $|X_n-X_n'|\le \delta_0$. Applying the above observations, 
$$
b_{k,n}^{-1} \E\big[|\eta_{k,n}^{(+,\ms{r})}(H_\ell) - \eta_{k+1,n}^{(-,\ms{r})}(H_\ell)|\one_{A_n}\big] \le \delta b_{k,n}^{-1} \P(A_n) + C^* b_{k,n}^{-1} \P \big( A_n\cap \big\{ |X_n-X_n'|>\delta_0 \big\} \big). 
$$
By \eqref{e:vague.Gkn1}, 
$$
\limsup_{n\to\infty}\delta b_{k,n}^{-1} \P(A_n) \le \delta\lim_{n\to\infty} b_{k,n}^{-1}\P(G_{k,n}^+(u_0)\ge 1)=\delta D_k e^{-u_0}. 
$$
As $\delta>0$ is arbitrary, it now  remains to show that for any $\delta_0 > 0$,  
$$
\lim_{n\to\infty} b_{k,n}^{-1} \P \big( A_n\cap \big\{ |X_n-X_n'|>\delta_0 \big\} \big)=0. 
$$
Let $\delta_0 > 0$ be fixed. For the proof, let us partition $[u_0,u_1)=\bigcup_{i=1}^m [v_i,v_{i+1})$ such that $|v_{i+1}-v_i|\le \delta_0$ for all $i\in\{ 1,\dots,m\}$. For $1\le i \le m$, set
\begin{align*}
G_{k,n}^+(i):= G_{k,n}^+(v_i) - G_{k,n}^+(v_{i+1}), \ \ \text{ and } \ \ G_{k+1,n}^-(i):= G_{k+1,n}^-(v_i) - G_{k+1,n}^-(v_{i+1}). 
\end{align*}
By construction, one can see that 
$$
A_n\cap \big\{ |X_n-X_n'|>\delta_0 \big\} \subset \bigcup_{i=1}^m \big\{  G_{k,n}^+(i) \neq G_{k+1,n}^- (i)\big\}. 
$$
Indeed, if $|X_n-X_n'|>\delta_0$ under $A_n$, then $X_n$ and $X_n'$ must fall into distinct subsets of the partition. Since $m$ is finite, it suffices  to demonstrate that for all $1 \leq i \leq m$,  as  $n \to \infty$,
$$
 b_{k,n}^{-1} \P\big(G_{k,n}^+(i) \neq G_{k+1,n}^-(i)\big) \to 0.
$$
This is however  a direct consequence of \eqref{e:pos.neg.diff1} and \eqref{e:pos.neg.diff2} by restricting the state space properly. 
\end{proof}

\bibliography{M0-critical}

\end{document}